\documentclass{amsart}
\usepackage{amsfonts,amscd}

\newtheorem{theorem}{Theorem}[section]
\newtheorem{lemma}[theorem]{Lemma}
\newtheorem{corollary}[theorem]{Corollary}
\newtheorem{proposition}[theorem]{Proposition}
\theoremstyle{remark}

\theoremstyle{definition}
\newtheorem{definition}[theorem]{Definition}

\numberwithin{equation}{section}
\makeatother

\DeclareMathOperator{\Rdb}{{\mathbb R}}

\DeclareMathOperator{\Ndb}{{\mathbb N}}

\begin{document}

\title[Characterizations of noncommutative $H^\infty$]{Characterizations
of noncommutative $H^\infty$}

\date{November 21, 2004.  Revised December 3, 2004.}

\author{David P. Blecher}
\address{Department of Mathematics, University of Houston, Houston, TX
77204-3008}
\email[David P. Blecher]{dblecher@math.uh.edu}
 \author{Louis E. Labuschagne}
\address{Department of Math, Applied Math, and Astronomy, P.O. Box
 392, 0003 UNISA, South Africa}
\email{labusle@unisa.ac.za}
\thanks{*Blecher is partially supported by a grant from
the National Science Foundation.
Labuschagne is  partially supported by the National Research
Foundation.}

\begin{abstract}
We transfer a large part of the circle of theorems characterizing
the generalization of classical $H^\infty$ known as `weak*
Dirichlet algebras', to Arveson's noncommutative setting of
subalgebras of finite von Neumann algebras.
\end{abstract}

\maketitle

\section{Introduction}

Around the early 1960's, it became apparent that a circle of famous
theorems about the classical $H^\infty$ space of bounded analytic
functions on the disk, could be generalized
to the setting of abstract function algebras.   In
particular, in \cite{Ho}, Hoffman showed that many of these
theorems were valid for function algebras satisfying the
{\em logmodular} condition
which he introduced in that paper.   In \cite{SW},
Srinivasan and Wang isolated a function algebra setting in which
the conclusions of many of these theorems, were each
equivalent
to logmodularity, and hence equivalent to each other.
  They called the algebras satisfying these
equivalent conditions {\em weak* Dirichlet algebras}.
Just a few years later, Arveson  introduced
his `finite maximal subdiagonal algebras' \cite{AIOA},
which we will consistently refer
to as {\em noncommutative $H^\infty$ algebras} in our paper,
for the sake of simplicity and brevity.
The setting, and we will fix this notation for the rest of
our paper, is a von Neumann algebra $M$
possessing a faithful normal tracial state $\tau$.
By a noncommutative $H^\infty$ algebra, we shall mean a
subalgebra of $M$ satisfying certain conditions
that we shall spell out momentarily.
These were intended to be
the noncommutative generalization of weak* Dirichlet algebras.
 One may then consider the possible noncommutative versions
of the famous
theorems about classical $H^\infty$, and ask which of
the conclusions of these theorems
are equivalent to the conditions defining the
noncommutative $H^\infty$ algebras.   This is the topic of
the present paper.

In \cite{BL}, we defined
a subalgebra $A$ of a $C^*$-algebra $B$
to be {\em logmodular} if every
strictly positive element $b \in B$
(that is, every selfadjoint $b$ for which there
exists an $\epsilon > 0$ with
$b \geq \epsilon 1$), is a uniform limit of
terms of the form $a^* a$ where $a \in A^{-1}$.
Here, $A^{-1}$ is the set of invertible elements of $A$.
We say that $A$ has {\em factorization}, if each strictly
positive $b \in B$ may be written as $a^* a$ for some $a \in A^{-1}$.
In \cite{BL}, we also defined
a {\em tracial subalgebra} of the algebra  $M$ above,
to be a weak* closed unital
subalgebra $A$ of $M$ for which there exists a projection $\Phi$
from $A$ onto ${\mathcal D} \overset{{\rm def}}{=} A \cap A^*$, such that
$\Phi$ is also a homomorphism, and $\tau = \tau \circ \Phi$ on $A$.
By Theorem 5.6 of \cite{BL},  $\Phi$ is
precisely the restriction to $A$ of the unique faithful normal
conditional expectation $\Psi$
from $M$ onto ${\mathcal D}$ such that $\tau = \tau \circ \Psi$.
Hence we may continue to write $\Psi$ as $\Phi$,
and we call this extension the
{\em conditional expectation} onto ${\mathcal D}$.

If ${\mathcal S}$ is a set, we will write
${\mathcal S}^*$ for the set $\{a :
 a^* \in {\mathcal S} \}$.
One may define  a noncommutative $H^\infty$ algebra to
be a tracial subalgebra  of $M$ for which
$A + A^*$ is weak* dense in $M$.
Arveson showed that noncommutative $H^\infty$ algebras have factorization,
and thus they are logmodular; in the present paper we will
prove the converse.
Also, Arveson formulated the classical Szeg\"o theorem,
and the related Jensen's inequality, in terms
of the  Fuglede-Kadison determinant.
This is a certain function
$\Delta : M \to [0,\infty)$, which
we will describe explicitly in Section 1.  Following Arveson, we
say that a tracial subalgebra $A$ satisfies
\begin{itemize}
\item {\em Jensen's inequality}, if $\Delta(\Phi(a)) \leq \Delta(a)$ for
all $a \in A$,
\item {\em Jensen's formula}, if $\Delta(\Phi(a)) = \Delta(a)$ for
all $a \in A^{-1}$,
\item {\em Szeg\"o's theorem}, if $\Delta(h) = \inf \{ \tau(h |a+d|^2)
: a \in A_0, d \in {\mathcal D}^{-1} , \Delta(d) \geq 1 \}$, for
all $h \in L^1(M)_+$.
\end{itemize}
Here and in the rest of our paper, $A_0 = A \cap {\rm Ker}(\Phi)$,
an ideal in $A$; and $L^1(M)_+$ is the positive part of the predual
of $M$.
We remark that our formulation of Szeg\"o's theorem may look different
from Arveson's formulation, but it is in fact equivalent
(see the Remark before Lemma \ref{Delgeq}).
Arveson had shown that for noncommutative $H^\infty$ algebras,
Jensen's inequality, Jensen's formula,
and his version of Szeg\"o's theorem,
are all equivalent.
Recently, the second author showed in \cite{LL3}
 that for noncommutative $H^\infty$
these three conditions are true (this was a major problem left
open in \cite{AIOA}).   For a general tracial subalgebra of $M$,
these conditions
are certainly not equivalent\footnote{One can see this
from the fact that ${\mathcal D}$ satisfies Jensen's inequality and
Jensen's formula, but not Szeg\"o's theorem in general.},
but at least Szeg\"o's theorem implies the others.
We show
that a tracial subalgebra satisfying Szeg\"o's theorem,
also satisfies the noncommutative variant of another of the
important
equivalent characterizations of weak* Dirichlet algebras,
namely:
\begin{equation} \label{cond}
{\rm If} \; g \in L^1(M)_+ \; {\rm and} \; \tau(fg) = \tau(f) \;
\textrm{for all} \; f \in
A, \; {\rm then} \; g = 1.
\end{equation}
This  condition may be rephrased as saying that there is
a unique normal state on $M$ extending $\tau_{\vert A}$.
Thus, we say that an algebra {\em has the unique normal state extension
property}
if it satisfies (\ref{cond}).
In the classical situation, (\ref{cond}) is exactly the key
condition underpinning the paper \cite{HR}.
By analogy with the work of
Hoffman and Srinivasan and Wang, one would expect to be able to
complete this circle; namely that any tracial subalgebra of $M$
having the unique normal state extension property,
is a noncommutative $H^\infty$ algebra.
We are able to show this under extra hypotheses,
and in particular
if $A$ is {\em $\tau$-maximal}, by which we mean that
$$A = \{ x \in M : \tau(
x A_0) = 0 \}.$$

 The main result of the present paper is the following:

\begin{theorem} \label{main}
 For a tracial subalgebra $A$ of $M$,
the following conditions
are equivalent:
\begin{itemize} \item [(a)]  $\overline{A + A^*}^{w*} = M$
(that is, $A$ is a noncommutative
$H^\infty$ algebra),
\item [(b)]   $A$ has factorization,
\item [(c)]  $A$ is logmodular,
\item [(d)]  $A$ satisfies Szeg\"o's theorem,
\item [(e)]  $A$ is $\tau$-maximal and has the
unique normal state extension property,
  \item [(f)]  $A$ has the unique normal state extension property, and
 $A + A^*$ is norm-dense in $L^2(M)$ (the latter space is
defined below).
\end{itemize}
\end{theorem}

At the end of
our paper, we briefly discuss characterizations of noncommutative
$H^\infty$ algebras in terms of an invariant subspace theorem,
or a Beurling-Nevanlinna factorization result.
These matters are essentially disjoint from
the rest of our paper, and deserve a more detailed
investigation at some point in the future.

We end this introduction with a few other
notational conventions.
We write $B^{-1}$ for the set of invertible elements of an algebra
$B$.  For a set ${\mathcal S}$, we write ${\mathcal S}_+$ for the set
$\{ x \in {\mathcal S} : x \geq 0 \}$,
where the symbol `$\geq$'
will usually denote the natural ordering in a $C^*$-algebra,
or in the predual of a von Neumann algebra, or more generally in the
noncommutative $L^p$ spaces.
We recall the definition of the latter:
For our (finite) von Neumann algebra $M$,
 we define
$\widetilde{M}$ to be
the set of closed densely defined operators $a$
affiliated to $M$.
If $1 \leq p < \infty$,  then $L^p(M,\tau) =
 \{a \in \widetilde{M} :  \tau(|a|^p) < \infty\}$,
equipped with the norm
$\|\cdot\|_p = \tau(|\cdot|^p)^{1/p}$.
 For brevity, we will in the following write $L^p$ or $L^p(M)$ for
$L^p(M,\tau)$, and we recall that for $M$ as above,
$L^p(M)$ may be defined to be the completion of $M$ in the
norm $\tau(\vert \cdot \vert^p)^{\frac{1}{p}}$.
The spaces $L^p$ are Banach spaces satisfying the usual duality
relations and H\" older inequalities \cite{FK,PX}. Now let $\Phi$
be a faithful normal conditional expectation from $M$ onto a von
Neumann subalgebra ${\mathcal D}$ satisfying $\tau \circ \Phi = \tau$.
The argument employed in e.g.\ Proposition 3.9 of \cite{MW}, then shows
that for each $1 \leq p < \infty$, the map $\Phi$ continuously
extends to a map which contractively maps $L^p(M)$ onto  $L^p({\mathcal D})$.

Note that $M^{-1} \cap M_+$ is precisely the set of
strictly positive elements of $M$.   We will use the fact
that in a finite von Neumann algebra, $x y = 1$ implies that
$y x = 1$.  We write $[{\mathcal S}]_p$ for the closure of a set
${\mathcal S}$ in $L^p(M)$.  If $A$ is a
noncommutative $H^\infty$ algebra, then
$[A]_p$ is often called the {\em
noncommutative Hardy space} $H^p$.

\section{The determinant}

For $a \in M$, the determinant
$\Delta(a)$ is
defined as follows: if $|a|$ is strictly positive,
then we define $\Delta(a) = \exp \tau(\log |a|)$.      Otherwise,
we define $\Delta(a) = \inf \,  \Delta(|a| + \epsilon 1)$,
the infimum taken over all scalars $\epsilon > 0$.

We will use several
basic properties of this determinant from \cite{AIOA}.
We will also need to extend the definition of this determinant to $L^1(M)$.
Namely, for any $a \in L^1(M)$, we set
 $\Delta(a) = \exp \tau(\log |a|)$ if $|a| \geq \epsilon 1$ for some
scalar $\epsilon > 0$; and otherwise,
we define $\Delta(a) = \inf \,  \Delta(|a| + \epsilon 1)$,
the infimum again taken over all scalars $\epsilon > 0$.
Clearly $\Delta(|a| + \epsilon 1)$ decreases as $\epsilon \to 0$,
with limit $\Delta(a)$.

We now explain why our definition of $\Delta(a)$ makes sense.  If $0 < \epsilon
< 1$, the function $\log t$ is bounded on $[\epsilon,1]$.
Also, $0 \leq \log t \leq t$ for $t \in [1,\infty)$.  So given
$h \in L^1(M)_+$ with $h \geq \epsilon 1$, it follows from the
Borel functional calculus for positive unbounded operators that
$(\log h) e_{[0,1]}$ is bounded, and that
$0 \leq (\log h) e_{[1,\infty)} \leq h e_{[1,\infty)}  \leq h.$  Here
$e_{[0,\lambda)}$ denotes the spectral resolution for $h$.   Thus
$(\log h) e_{[0,1]}$ and $(\log h) e_{[1,\infty)}$ belong to
$L^1(M)$.  Hence $\log h \in L^1(M)$ in this case.

We will need the following
variant of some formulae from \cite{AIOA}:

\begin{proposition} \label{Arvf}  {\rm (Cf.\ 4.3.2 in \cite{AIOA})} \
For any $h \in L^1(M)_+$, we have
\begin{equation} \label{Eq1}
\Delta(h) \; = \; \inf \{
\tau(hb) : b \in M_+ \cap M^{-1} , \, \Delta(b) \geq 1 \} .
\end{equation}  Also, this infimum is attained on the
von Neumann algebra generated by $h$ (see e.g.\ \cite[p.\ 349]{KR}), which is
a commutative subalgebra of $M$.
For any $h \in L^1(M)$, we have
\begin{equation} \label{Eq2}
\Delta(h) \; = \; \inf \{ \tau(|hb|) : b \in M_+ \cap M^{-1} , \,
 \Delta(b) \geq 1 \} .
\end{equation}
\end{proposition}

\begin{proof}  First assume that $h \geq \epsilon 1$ for some
$\epsilon > 0$.  As in the calculation above,
$$0 \leq (\log h) e_{[n,\infty)} \leq h e_{[n,\infty)}
\; , \qquad n \in \Ndb .$$
Since $e_{[n,\infty)} \to 0$ strongly, and
since $h \in L^1(M)_+$, we have $\lim_n \tau(h e_{[n,\infty)}) = 0$.
Thus $\lim_n \tau((\log h) e_{[n,\infty)})  = 0$.
Hence $\tau(\log h) = \lim_n \tau((\log h) e_{[0,n)})$.
  That is, $$\Delta(h) \; = \; \exp \tau(\log h)   \; = \;  \lim_n
 \exp \tau((\log h) e_{[0,n)}) .$$

Given any $b \in M_+$, we have
$$\tau(h e_{[0,n)} b)  \; = \; \tau(b^{\frac{1}{2}} h e_{[0,n)}
b^{\frac{1}{2}})  \; \leq \; \tau(b^{\frac{1}{2}} h b^{\frac{1}{2}})
  \; = \; \tau(h b) .$$
Combining this fact with \cite[4.3.2]{AIOA}, we have that
\begin{align*}
\Delta(he_{[0,n)}) & \; = \;
\inf \{ \tau(he_{[0,n)} b) : b \in M_+ \cap M^{-1} , \, \Delta(b) \geq 1 \} \\
& \; \leq \;
 \inf \{ \tau(h b) : b \in M_+ \cap M^{-1} , \, \Delta(b) \geq 1 \} .
\end{align*}
Now let $n \to \infty$ to see that
$$\Delta(h) \; \leq \inf \{ \tau(h b) : b \in M_+ \cap M^{-1} ,
 \, \Delta(b) \geq 1 \} .$$
To see that the infimum is precisely $\Delta(h)$, and that it
is attained on the minimal commutative
von Neumann algebra $M_0$ generated by $h$, it suffices to find,
for each scalar $\delta > 0$, an element $b_\delta$ in
this von Neumann algebra, with  $\Delta(b_\delta) \geq 1$,
and $\inf \{ \tau(h b_\delta) : \delta > 0 \} = \Delta(h)$.
Since $h  \geq \epsilon 1$, we have that $h^{-1}$ is bounded, and is
in $M_+$.  Also $\infty > \Delta(h) > 0$, since $\log h \in L^1(M)$.
Set $b_\delta = \Delta(h) (h^{-1} + \delta 1)$
for each $\delta > 0$.    Then $b_\delta \in M_+ \cap M^{-1}$,
and indeed $b_\delta \in M_0$.  Clearly
$$\inf \{ \tau(h b_\delta) : \delta > 0 \}
 \; = \; \Delta(h) \inf \{ \tau(1 + \delta h)  : \delta > 0 \}  \; =
\; \Delta(h) \; .$$
To see that $\Delta(b_\delta) \geq 1$, first note that
both $\log h$ and $\log(h^{-1} + \delta 1)$ are in
$L^1(M)$, by virtue of the fact that $h \geq \epsilon 1$ and
$h^{-1} + \delta 1 \geq  \delta 1$.
Then \begin{align*}
\Delta(b_\delta) & = \Delta(h) \Delta(h^{-1} + \delta 1)
= (\exp \tau(\log(h)))(\exp \tau(\log(h^{-1} + \delta 1))) \\
&  = \exp \tau(\log(h) + \log(h^{-1} + \delta 1)) .
\end{align*}
 By the Borel functional
calculus for unbounded selfadjoint operators,
$$0 \; \leq \; \log(1 + \delta h) \; = \;
 \log(h(h^{-1} + \delta 1))
\; = \; \log(h) + \log(h^{-1} + \delta 1) , $$
the sum and product here being the `strong' ones (that is,
we are taking closures of the operators involved).  Hence,
$$\Delta(b_\delta) = \exp \tau(\log(1 + \delta h)) \geq 1 .$$

Next, if $h$ is not $\geq \epsilon 1$ for any $\epsilon > 0$,
then by the above we have
$$\Delta(h) = \inf_{\epsilon > 0} \Delta(h + \epsilon 1)
= \inf_{\epsilon > 0} \inf \{ \tau((h + \epsilon 1)b) : b \in M_+ \cap M^{-1} , \, \Delta(b) \geq 1 \} .$$
Interchanging the infimums gives (\ref{Eq1}).
The infimum is again achieved on the von Neumann algebra
generated by $h$, since for any $\epsilon > 0$ this is the same as the
 von Neumann algebra generated by $h  + \epsilon 1$.

To obtain (\ref{Eq2}), first note that
we may assume that $h \geq 0$.  This is because if the result held in
the latter case, then in the general case,
\begin{align*}
\Delta(h)  & = \Delta(|h|) = \inf \{ \tau(||h|b|) : b \in M_+ \cap M^{-1}
\, , \, \Delta(b) \geq 1 \} \\
& = \inf \{ \tau(|h b| : b \in M_+ \cap M^{-1} \, , \, \Delta(b) \geq 1 \} ,
 \end{align*} since $|h b|^2 = b h^* h b = ||h|b|^2$.

For $h \geq 0$, it then follows from the above that
\begin{align*} \Delta(h) & = \inf \{ \tau(hb) : b \in M_+ \cap M^{-1} ,
\, \Delta(b) \geq 1 \} \\ & \leq
\inf \{ \tau(|hb|) :  b \in M_+ \cap M^{-1} ,
\, \Delta(b) \geq 1 \}.  \end{align*}
To see that the last $\leq$ is an equality, we
argue as in the corresponding
part of the proof above.  Indeed, if $\epsilon, \delta > 0$ then
\begin{align*} \Delta(h) & \leq \tau(|h \Delta(h + \epsilon 1)
((h + \epsilon 1)^{-1} + \delta 1)|)  =
\tau(h \Delta(h + \epsilon 1)
((h + \epsilon 1)^{-1} + \delta 1)) \\
& \leq \tau((h + \epsilon 1)  \Delta(h + \epsilon 1)
((h + \epsilon 1)^{-1} + \delta 1))  \, \longrightarrow \,
\Delta(h)
\end{align*}
as $\delta, \epsilon \to 0$ (similarly to the first part of the
proof).
This establishes (\ref{Eq2}).
   \end{proof}

{\bf Remark.}  In \cite{AIOA},     Arveson defined the quantity
$\Delta(\rho)$, for a normal state $\rho$ of $M$, to be
$\inf \{ \rho(b) : b \in M_+ \cap M^{-1} ,
\, \Delta(b) \geq 1 \}$.  However, normal states $\rho$ are
in bijective correspondence with the norm-one elements
$h$ of $L^1(M)_+$,  via the
map $h \mapsto \rho_h = \tau(h \cdot)$.
We may therefore rephrase (\ref{Eq1})
as the statement that
$\Delta(h) = \Delta(\rho_h)$ for all such $h$.

The fact that our formulation of
the noncommutative  Szeg\"o  theorem in
Section 1, is in line with that of Arveson, follows
from this, and from the following Claim:
\begin{align*} \inf \{ \tau(h |a+d|^2) :  & a \in A_0, d \in {\mathcal D}^{-1} ,
\Delta(d) \geq 1 \}  \\
&  = \; \inf \{ \tau(h |a+d|^2) :  a \in A_0, d
\in {\mathcal D} , \Delta(d) \geq 1 \} ,
\end{align*}
 for any $h \in L^1(M)_+$.
 Indeed, clearly the left side dominates
the right.  To prove the converse inequality, first note that
the infimum on the right side is unchanged if one insists also that
$d \geq 0$.  For if $d
\in {\mathcal D}$, with $\Delta(d) \geq 1$, let $d = v|d|$ be the polar
decomposition of $d$ in ${\mathcal D}$.  Set $\widetilde{a}
= v^*a \in A_0$. We have $\Delta(|d|) = \Delta(d) \geq 1$.
 Also $|a + d|^2 \geq |\widetilde{a} +
|d||^2$, so that $\tau(h|a + d|^2) \geq \tau(h|\widetilde{a} +
|d||^2)$.

If $d \in {\mathcal D}_+$, then $d + \frac{1}{n} \in
{\mathcal D}^{-1}$, and $\Delta(d + \frac{1}{n}) \geq \Delta(d) \geq 1$,
and
$$\tau(h|a + d|^2) \; = \;  \lim_n \tau(h|a + d + \frac{1}{n}|^2)
\; \geq \; \inf_n \tau(h|a + d + \frac{1}{n}|^2) ,$$
 which dominates the left side of the equation in the Claim.

\begin{lemma}   \label{Delgeq}  Let $h \in L^1(M)_{sa}$.
If for some $\delta > 0$ we have
$$\Delta(1 - th) \geq 1 \; , \qquad t \in (-\delta,\delta),$$
then $h = 0$.
\end{lemma}

\begin{proof}  Let $M_0$ be the von Neumann algebra generated by
$h$ (see e.g.\ \cite[p.\ 349]{KR}), which is
a commutative subalgebra of $M$.
Let $\psi = \tau_{\vert M_0}$.
Since $\psi$ is a faithful normal state on $M_0$, it is a simple
consequence of the Riesz representation theorem applied to $\psi$, that
$M_0 \cong L^\infty(\Omega,\mu_\tau)$ $*$-algebraically, for a
measure space $\Omega$ and a Radon probability
measure $\mu_\tau$. (This is a simpler case of the proof of 1.18.1
 in \cite{Sak}).
We have $\tau(x) = \int_\Omega \, x \, d \mu_\tau$ for any
$x \in M_0$, where we are abusing notation by writing $x$ for
the corresponding element of $L^\infty(\Omega,\mu_\tau)$ too.
Then $L^1(M_0)$, which is the completion of $M_0$ in the
norm $\tau(|\cdot|)$, is isometric to
$L^1(\Omega,\mu_\tau)$, the completion of $L^\infty(\Omega,\mu_\tau)$
in the $L^1$-norm.   Clearly,
$L^1(M_0) \subset L^1(M)$ isometrically,
and the canonical extension of $\tau$ to $L^1(M)$,
agrees with the canonical extension of $\psi$ to
$L^1(M_0)$.  In particular, for $h \in L^1(M_0)$,
$\tau(h) =  \int_\Omega \, h  \, d \mu_\tau$,
and $\tau(|h|)$ is the $L^1(\mu_\tau)$-norm of $h$.

Since we may compute the quantity $\Delta(1-th)$ with
respect to $M_0$, it follows from the last paragraph, and
Lebesgue's monotone convergence theorem,  that
 $$\Delta(1-th) \; = \; \inf_{\epsilon > 0}
 \{ \exp \int_\Omega \, (\log |1-th| +
\epsilon 1) \, d
\mu_\tau  \}   \; = \;   \exp \int_\Omega \, \log |1-th|
\, d
\mu_\tau .$$
Thus  the result we want follows from Hoffman's  lemma
 \cite[Lemma 6.6]{Ho}.
 \end{proof}

\section{Consequences of logmodularity}

\begin{proposition} \label{Jenlog}  For a tracial
subalgebra, logmodularity implies Jensen's formula.
\end{proposition}

\begin{proof}  This is a modification of the proof
of the main result of \cite{LL3}.  We simply indicate
the parts of the proof which need adjusting, beginning
at the inductive step.  We assume that
$$\tau(|a|^{\frac{1}{2^k}}) \; \geq \; \tau(|\Phi(a)|^{\frac{1}{2^k}})
\; , \qquad a \in A^{-1} ,$$
for an integer $k$, and we want to prove the same
inequality with $k$ replaced by $k+1$.  Fix
$a \in A^{-1}$, and inductively
define $(x_n) \subset M_+$ by
$$x_1 = |a|^{\frac{1}{2^k}} \; , \; x_{n+1} = \frac{1}{2}(x_n +
|a|^{\frac{1}{2^k}} x_{n}^{-1}) .$$
For each $n \in \Ndb$, select $(z_m^{(n)}) \subset A^{-1}$ with
$$\lim_m \, |z_m^{(n)}| \; = \; x_n^{2^k} \; ,
 \qquad n \in \Ndb.$$
For $q = \frac{1}{2^k}$, we have, as in \cite{LL3}, that
\begin{equation} \label{one} \tau(|a (z_m^{(n)})^{-1}|^q)
\; = \;  \tau((|a||z_m^{(n)}|^{-2} |a|)^{\frac{q}{2}})
 \; \qquad m \in \Ndb.
\end{equation}
By the inductive hypothesis, and H\"older's inequality,
we have
\begin{align*} \frac{1}{2} \tau(|z_m^{(n)}|^{q} +
|a (z_m^{(n)})^{-1}|^{q}) &
\geq \frac{1}{2} \tau(|\Phi(z_m^{(n)})|^{q} +
|\Phi(a) \Phi(z_m^{(n)})^{-1}|^{q}) \\
& \geq \frac{1}{2} \Big\{ \tau(|\Phi(z_m^{(n)})|^{q}) +
\frac{\tau(|\Phi(a)|^{\frac{1}{2^{k+1}}})^2}{\tau(
|\Phi(z_m^{(n)})|^{q})} \Big\} \\
& \geq \tau(|\Phi(a)|^{\frac{1}{2^{k+1}}}) .
\end{align*}
Using the last inequality, and (\ref{one}), we obtain
$$\tau(|\Phi(a)|^{\frac{1}{2^{k+1}}}) \leq
\frac{1}{2} \tau(|z_m^{(n)}|^{\frac{1}{2^k}} +
|a (z_m^{(n)})^{-1}|^{\frac{1}{2^k}})
= \frac{1}{2} \tau(|z_m^{(n)}|^{\frac{1}{2^k}} +
                                 (|a| |z_m^{(n)}|^{-2} |a|)^{\frac{1}{2^{k+1}}}).
$$
The left side of this inequality does not depend on $m$.
Letting $m \to \infty$ gives
$$\tau(|\Phi(a)|^{\frac{1}{2^{k+1}}}) \leq
\frac{1}{2} \tau(x_n + x_n^{-1} |a|^{\frac{1}{2^k}}) =
\tau(x_{n+1}) .$$
Since the $x_{n+1}$'s decrease monotonically to $|a|^{\frac{1}{2^{k+1}}}$,
 we have
$\tau(|\Phi(a)|^{\frac{1}{2^{k+1}}}) \leq \tau(|a|^{\frac{1}{2^{k+1}}})$,
 as required.
 \end{proof}

\begin{proposition} \label{islo}  Let $A$ be a logmodular tracial
subalgebra of $M$.  For any $h \in L^1(M)$, we have
$$\Delta(h) \; = \; \inf \{ \tau(|ha|) : a \in A^{-1} ,
\Delta(a) \geq 1 \} .$$
\end{proposition}

\begin{proof}    Since $|ha|^2 = a^* |h|^2 a = ||h| a|^2$,
we have
$$ \inf \{ \tau(|ha|) : a \in A^{-1} ,
\Delta(a) \geq 1 \} \; = \;  \inf \{ \tau(||h| a|) : a \in A^{-1} ,
\Delta(a) \geq 1 \} .$$
Since $\Delta(h) = \Delta(|h|)$, it suffices henceforth to
assume that $h \geq 0$.

Given $a \in A$, we have that
\begin{equation} \label{sta} \tau(|ha|)
= \tau(|a^* h|) = \tau(||a^*| h|) = \tau(|h |a^*||) .
\end{equation}
If $a \in A^{-1} \subset M^{-1}$, then $a^*$ and $|a^*|$
are in $M^{-1}$.   Since $\Delta(a) = \Delta(a^*)
= \Delta(|a^*|)$, by \cite[4.3.1]{AIOA}, it follows from (\ref{Eq2})
 that
\begin{align*} \inf \{ \tau(|ha|) : a \in A^{-1} ,
\Delta(a) \geq 1 \} & = \inf \{ \tau(|h |a^*||) :
 a \in A^{-1} ,
\Delta(a) \geq 1 \} \\
& \geq  \inf \{ \tau(|h b| :
b \in M_+ \cap M^{-1} ,
\Delta(b) \geq 1 \} \\
& = \Delta(h).
\end{align*}

Since $A$ is
logmodular, given any $b_0 \in M_+ \cap M^{-1}$ with
$\Delta(b_0) \geq 1$, we may select
$\{ a_n \} \subset A^{-1}$ so that $|a_n^*| \to b_0$
uniformly.  By \cite[4.3.1 (iv)]{AIOA}, we deduce that
$$\Delta(a_n) = \Delta(a_n^*) = \Delta(|a_n^*|) \; \longrightarrow \;
\Delta(b_0) .$$
Letting $\tilde{a}_n = \frac{\Delta(b_0)}{\Delta(a_n)} \, a_n \in
A^{-1}$, we have $|\tilde{a}_n^*| =
\frac{\Delta(b_0)}{\Delta(a_n)} |a_n^*| \to b_0$; and
$\Delta(\tilde{a}_n) = \Delta(b_0) \geq 1$.   Thus we clearly
have, using also (\ref{sta}), that
$$\tau(|h b_0|) = \lim_n \tau(|h |\tilde{a}_n^*||)
=  \tau(|h \tilde{a}_n|)
\geq \inf \{ \tau(|ha|) : a \in A^{-1} , \Delta(a) \geq 1 \}.$$
Since $b_0 \in M_+ \cap M^{-1}$ was arbitrary, the above
combined with the earlier inequality in this proof, gives
$$\Delta(h) \; = \;  \inf \{ \tau(|h b| :
b \in M_+ \cap M^{-1} ,
\Delta(b) \geq 1 \}  \; = \;
\inf \{ \tau(|ha|) : a \in A^{-1} ,
\Delta(a) \geq 1 \} ,$$
which is the desired equality.
     \end{proof}

\begin{lemma} \label{logimdel}  Let $A$ be a
logmodular tracial subalgebra of $M$.
If $h \in L^1(M)$ with $\tau(ha) = 0$ for
every $a \in A$, then $\Delta(1-h) \geq 1$.
\end{lemma}

\begin{proof}
Suppose that $\tau(ha) = 0$ for
every $a \in A$.   We continue to write
$\Phi$ for the canonical `extension by continuity' of
$\Phi$ to a map from $L^1(M)$ to $L^1({\mathcal D})$
(see e.g.\ \cite{Sai} or 3.9 in \cite{MW}).
By routine approximation arguments, it is easy to
see that this extension
is still a contractive `conditional expectation':
$\Phi(h d) = \Phi(h) d$ for $h  \in L^1(M), d \in {\mathcal D}$.
Similarly, $\tau \circ \Phi = \tau$ on $L^1(M)$.
Using these facts, given $a_0 \in A, d \in {\mathcal D}$,
we have
\begin{align*} \tau(\Phi((1-h) a_0) d) & = \tau(\Phi((1-h) a_0 d))
=  \tau((1-h) a_0 d) \\
& = \tau(a_0 d) = \tau(\Phi(a_0 d)) = \tau(\Phi(a_0) d) .
\end{align*}
This implies that $\Phi((1-h) a_0) = \Phi(a_0)$.
Next, let $a_0 \in A^{-1}$ be given.
Since $\Phi$ is contractive on $L^1$, we have, by (\ref{Eq2}) and
 the above, that
\begin{align*} \tau(|(1-h) a_0|) \; & \geq \;
\tau(|\Phi((1-h) a_0)|) = \tau(|\Phi(a_0)|) \\
& \geq \; \inf \{ \tau(|\Phi(a_0) b|)
 : b \in M_+ \cap M^{-1} , \Delta(b) \geq 1 \} \\
 & = \;  \Delta(\Phi(a_0)) .
\end{align*}
By Jensen's formula, the latter
quantity equals $\Delta(a_0)$.
By Proposition \ref{islo}, we conclude that $\Delta(1-h) \geq 1$.
\end{proof}

\begin{corollary} \label{logeqnon}  For a
tracial subalgebra $A$ of $M$, the following are equivalent:
\begin{itemize} \item [(a)]  $\overline{A + A^*}^{w*} = M$
(that is, $A$ is a noncommutative
$H^\infty$ algebra),
\item [(b)]   $A$ has factorization,
\item [(c)]  $A$ is logmodular.
\end{itemize} \end{corollary}

\begin{proof}
That (a) $\Rightarrow$ (b) was
proved in \cite[4.2.1]{AIOA}.  That (b) $\Rightarrow$ (c) is
evident from the definition \cite{BL}. Let $h \in L^1(M)$ be given.
 We will write $h \in (A + A^*)_\perp$ for the claim that $\langle h,
 a\rangle = \tau(ha^*) = 0$ for each $a \in A + A^*$   If  $A$ is
logmodular, to show that $A + A^*$ is weak* dense in $M$, it suffices
 to show that if $h \in (A + A^*)_\perp$
then $h = 0$.  Since $A + A^*$ is
selfadjoint, it is easy to see that
$h \in (A + A^*)_\perp$ if and only if
$h + h^* \in (A + A^*)_\perp$ and $\frac{i}{2} (h-h^*)
 \in (A + A^*)_\perp$.
We may therefore assume that $h$ is selfadjoint.
By Lemma \ref{logimdel}, $\Delta(1-th) \geq 1$ for every $t \in
\Rdb$.   By Lemma \ref{Delgeq}, $h = 0$.
\end{proof}

It follows from Corollary \ref{logeqnon} and \cite{LL3},
that a logmodular tracial subalgebra $A$ satisfies
Jensen's inequality and
Szeg\"o's theorem.  In fact, one can prove this directly.
Indeed a more precise result is stated next, for which we
will need the following definitions:
$${\mathcal S}_1 = \{ b : b \in M^{-1} \cap M_+ , \,
\Delta(b) \geq 1 \} ,$$
$$ {\mathcal S}_2 = \{ a^* a : a \in A , \,
\Delta(\Phi(a)) \geq 1 \} , $$
$$ {\mathcal S}_3 = \{ a^* a : a \in A^{-1} , \,
\Delta(a) \geq 1 \} .$$

\begin{proposition} \label{mlm} \begin{itemize}
\item [(a)]   For any tracial subalgebra,
Szeg\"o's theorem implies
Jensen's inequality,
and Jensen's inequality implies
 Jensen's formula.  If $A$ is logmodular,
then Szeg\"o's theorem holds.

\item [(b)]   For any tracial subalgebra $A$,
${\mathcal S}_3 \subset {\mathcal S}_1$; and
$A$ is logmodular if and only if $\overline{{\mathcal S}_1}
= \overline{{\mathcal S}_3}$.
\end{itemize}   \end{proposition}

\begin{proof}
(b) \ It is clear that ${\mathcal S}_3 \subset {\mathcal S}_1$.
If $A$ is logmodular, and $b \in {\mathcal S}_1$,
there is a sequence $(c_m) \subset A^{-1}$ with $c_m^* c_m
\to b$.  By \cite[4.3.1]{AIOA}, $\Delta(c_m)^2
= \Delta(c_m^* c_m) \to \Delta(b)$.   Letting
$d_m = \frac{\sqrt{\Delta(b)}}{\Delta(c_m)} \; c_m$,
we have $d_m^* d_m \to b$, and
$\Delta(d_m) = \sqrt{\Delta(b)} \geq 1$.
Thus $b \in \overline{{\mathcal S}_3}$.  Thus
$\overline{{\mathcal S}_1}
= \overline{{\mathcal S}_3}$.

 Conversely, if $\overline{{\mathcal S}_1}
= \overline{{\mathcal S}_3}$, then a simple
approximation argument shows that $A$ is logmodular.
Namely, note that if $b \in M^{-1} \cap M_+$, we have
 $\Delta(b) \geq \epsilon\Delta(1) > 0$, where $\epsilon > 0$
 is chosen so that $b \geq \epsilon 1$.
Then by scaling, we may assume that $\Delta(b) = 1$.
Then $b \in {\mathcal S}_1 \subset  \overline{{\mathcal S}_3}$,
giving the desired approximations.  This
completes the proof of (b).

(a) \ Our proof is based on \cite[4.4.3]{AIOA}.
The fact  that
Szeg\"o's theorem
implies Jensen's inequality, follows from a slight modification of
the argument
Arveson used in \cite[4.4.3]{AIOA} to show that Szeg\"o's theorem implies
Jensen's formula.
 The main change is that we begin
the proof by choosing $a \in A$ (as opposed
to $A^{-1}$).   We obtain the inequality
$$\inf_{T} \; \rho(|D + T|^2) \; \geq \; \phi(|D \Phi(a)|^2),$$
in line 2 of p.\ 613 of \cite{AIOA}, instead of an
equality.  Following the argument until line 6 of
that page, we obtain $\Delta(\Phi(a)) \leq \Delta(a)$ as required.

To see that Jensen's inequality
implies Jensen's formula, is the calculation on lines 8 and
9 of p.\ 613 of \cite{AIOA}.

As noted on p.\ 612  of \cite{AIOA}, Arveson's
Szeg\"o theorem is
$$\inf \{ \tau(h b) : b \in {\mathcal S}_1 \} \; = \;
\inf \{ \tau(h b) : b \in {\mathcal S}_2 \} \; , \qquad
h \in L^1(M)_+ .$$
We will show that if $A$ is logmodular, then
$$\overline{{\mathcal S}_1}    \; = \;
 \overline{{\mathcal S}_2} \; = \;
 \overline{{\mathcal S}_3} .$$
This, together with the fact that
$\tau(h \cdot) $ is norm-continuous on $M$,
will complete the proof of (a).

If $A$ is logmodular, we have
${\mathcal S}_3 \subset {\mathcal S}_2$
by Proposition \ref{Jenlog}.
In the light of
 (b), to complete the proof we need only show that
${\mathcal S}_2 \subset \overline{{\mathcal S}_3}$.
To this end, let $a \in A$ with $\Delta(\Phi(a)) \geq 1$.
For fixed $n \in \Ndb$,  by logmodularity there
exists a sequence $(c_m) \subset A^{-1}$ with
$a^* a + \frac{1}{n} 1 = \lim_m \, c_m^* c_m$.
Thus for any $\epsilon > 0$ there exists an $N_\epsilon$
such that $a^* a + \frac{1}{n} 1  \leq (1+\epsilon) c_m^* c_m$
for all $m \geq N_\epsilon$.
By passing to a subsequence, we can assume that
\begin{equation} \label{aneq2}
a^* a + \frac{1}{n} 1  \; \leq \; (1 + \frac{1}{m}) c_m^* c_m \; ,
\qquad m \in \Ndb .
 \end{equation}
Since $\Phi$ is completely positive, we may use the
Kadison-Schwarz inequality (i.e.\ $\Phi(x)^* \Phi(x) \leq
\Phi(x^* x))$, to see that
\begin{equation} \label{aneq1}
\Phi(c_m^{-1})^* \Phi(a)^*
\Phi(a) \Phi(c_m^{-1})  = \Phi( a c_m^{-1})^* \Phi( a c_m^{-1}) \leq
\Phi((c_m^{-1})^* (a^* a + \frac{1}{n} 1) c_m^{-1}) .
\end{equation}
Combining (\ref{aneq1}) and  (\ref{aneq2}), we have
$$\Phi(c_m^{-1})^* \Phi(a)^*
\Phi(a) \Phi(c_m^{-1}) \leq \Phi((1 + \frac{1}{m}) 1) = (1 + \frac{1}{m}) 1 .$$
Left and right multiplying by $\Phi(c_m)^*$ and $\Phi(c_m)$ respectively,
we have
\begin{equation} \label{aneq3}
\Phi(a)^* \Phi(a) \; \leq \;  (1 + \frac{1}{m}) \Phi(c_m)^* \Phi(c_m)  \; ,
\qquad m \in \Ndb .                 \end{equation}
It follows from Proposition \ref{Jenlog} and  \cite[4.3.1]{AIOA},
that
$$ \Delta(c_m^* c_m) = \Delta(c_m)^2 =  \Delta(\Phi(c_m))^2
=  \Delta(\Phi(c_m)^* \Phi(c_m)) .$$
From this and (\ref{aneq3}),
and facts in \cite[4.3.1]{AIOA},
we deduce that
$$\Delta(a^* a + \frac{1}{n} 1) \; = \; \lim_m \Delta(c_m^* c_m)
\; \geq \;  \lim_m \, \frac{\Delta(\Phi(a)^* \Phi(a))}{1+ m^{-1}}
 \; = \;   \Delta(\Phi(a))^2  \geq 1 .$$
Thus $a^* a + \frac{1}{n} \in {\mathcal S}_1
\subset \overline{{\mathcal S}_3}$.  Taking the limit
over $n$, we find $a^* a \in \overline{{\mathcal S}_3}$
as desired.
\end{proof}

\section{The unique normal state extension property}

In  \cite{Lum}, Lumer showed the importance of
the `uniqueness of representing measure'
criterion to the generalized $H^p$ theory.
Shortly thereafter, Hoffman and Rossi
showed that in the setting considered in
\cite{SW}, condition (\ref{cond})
 characterized weak* Dirichlet algebras.
Note that
in their setting,
$\tau_{|A} = \Phi_{|A}$ is multiplicative.
Lumer, on the other hand, required
there to be a unique state on $L^\infty$
(that is, a unique probability measure on
the maximal ideal space of $L^\infty$)
extending $\tau_{|A}$.

In our noncommutative context, there are several
conditions, besides  (\ref{cond}),
which present themselves as generalizations
of the `unique extension' properties mentioned in the
last paragraph.  For example,  one could consider
the condition that there be a unique completely
positive extension of $\Phi_{|A}$.
 A stronger condition yet, is that every
completely contractive representation\footnote{See e.g.\ \cite{BLM}
for the definition of these terms, if needed.}
of $A$ has a unique completely positive extension to $M$.
It is known that noncommutative $H^\infty$
algebras satisfy this latter condition \cite{BL}.
Although this condition does have some interesting
consequences, we have not yet been able to
connect it  convincingly to other properties
considered in this paper.  Thus in this section,
we focus on the condition (\ref{cond}).

\begin{lemma}  \label{gcon}   For a tracial subalgebra of $M$,
the unique normal state extension property {\rm (\ref{cond})}
 is equivalent to:
\begin{equation} \label{con}
{\rm If} \; g \in L^1(M)_+ \; {\rm and} \; \tau(fg) = 0 \;
\textrm{for all} \; f \in
A_0, \; {\rm then} \; g \in L^1({\mathcal D}) .
\end{equation}
\end{lemma}

\begin{proof}   Suppose that  (\ref{con}) holds.
If  $h \in L^1(M)_+$ with $\tau(h a) = \tau(a)$ for all $a \in A$,
then taking $a \in A_0$, we conclude that $h \in L^1({\mathcal D})$.
We also have $\tau(h d)  = \tau(d)$ for all $d \in {\mathcal D}$,
which forces $h = 1$.

Conversely, suppose that  (\ref{cond}) holds,
and that we are given an $h \in L^1(M)_+$,
 such that $\tau(h a) = 0$ for all $a \in A_0$.
We may suppose that $h \geq 1$, by replacing $h$ with
$h + 1$ if necessary.  Then also $\Phi(h) \geq 1$.
 If $a \in A$, we then have
$$\tau(h \Phi(a)) = \tau(\Phi(h \Phi(a))) =
\tau(\Phi(h) \Phi(a)) = \tau(\Phi(\Phi(h) a))
= \tau(\Phi(h)a) .$$
In the last line we have used several properties of $\Phi$ which
are obvious for $\Phi$ considered as a map on $M$, and
which are easily verified for the extension of $\Phi$ to
$L^1(M)$.   Hence
$$\tau(ha) = \tau(h \Phi(a)) + \tau(h (a - \Phi(a))) = \tau(h \Phi(a))
= \tau(\Phi(h)a) .$$
Since $\Phi(h) \geq 1$, $\Phi(h)^{-1} \in {\mathcal D}_+$.
We have $\Phi(h)^{-\frac{1}{2}} a \Phi(h)^{-\frac{1}{2}} \in A$
for every $a  \in A$.  Setting $\tilde{h} = \Phi(h)^{-\frac{1}{2}}
h \Phi(h)^{-\frac{1}{2}}$, we have $\tilde{h} \in L^1(M)_+$,
and by the previous centered equation,
$$\tau(\tilde{h} a) = \tau(h (\Phi(h)^{-\frac{1}{2}} a \Phi(h)^{-\frac{1}{2}}))
= \tau(\Phi(h) \Phi(h)^{-\frac{1}{2}} a \Phi(h)^{-\frac{1}{2}})
= \tau(a) ,$$
   for any $a  \in A$.   By (\ref{cond}),
$\tilde{h} = 1$, so that $h = \Phi(h) \in L^1({\mathcal D})$.
\end{proof}

\begin{theorem}  \label{gcond}  Let $A$ be a
tracial subalgebra of $M$ which satisfies
Szeg\"o's theorem.  Then $A$ has the
unique normal state extension property {\rm (\ref{cond})}.
\end{theorem}

\begin{proof}  Suppose that we are given an $h \in L^1(M)_+$,
 such that $\tau(h a) = \tau(a)$ for all $a \in A$.
Then $\tau(h a) = 0$ for all $a \in A_0$,
and hence also for all $a \in A_0^*$, since
$\overline{\tau(h a^*)} = \tau(ha)$.
If $a  \in A_0$, and $d \in {\mathcal D}$, then
$$\tau(h |a+d|^2) =
\tau(h |a|^2 + h d^* a + h a^* d + h |d|^2)
= \tau(h^\frac{1}{2} |a|^2 h^\frac{1}{2} + |d|^2)
\geq \tau(|d|^2) .$$
Appealing to Szeg\"o's theorem, we deduce that
$$\Delta(h) \; = \; \inf \{ \tau(|d|^2) : d \in
{\mathcal D}^{-1}, \Delta(d) \geq 1 \} .$$
By (\ref{Eq1}), and by \cite[4.3.1]{AIOA},
we have
$$\tau(|d|^2) \geq \Delta(|d|^2) = \Delta(|d|)^2 = \Delta(d)^2 .$$
It follows that
$$\Delta(h) \; = \; \inf \{ \tau(|d|^2) : d \in
{\mathcal D}^{-1}, \Delta(d) \geq 1 \}   \; = \;  1  .$$
By hypothesis, we also have $\tau(h) = \tau(1) = 1$.
We now consider the von Neumann algebra $M_0$ generated
by $h$.  With notations as in the proof of
Lemma \ref{Delgeq}, we have
$\int_\Omega \; h \, d \mu_\tau =
1 = \exp(\int_\Omega \;  \log h  \, d \mu_\tau)$.
It is an elementary exercise in real analysis to show that
this forces $h = 1$.  (Letting $k = \log h$, we have $\int \, k = 0$ and
$\int \, e^k  = 1$.  If $r = e^k - k - 1$ then $r$ is a nonnegative
function, but $\int \, r = 0$.   Thus $r = 0$, forcing $k = 0$.)
 \end{proof}

\begin{definition} \label{el2den}
We say that a tracial subalgebra $A$ of
$M$ satisfies {\em $L^2$-density}, if
$A + A^*$ is dense in $L^2(M)$ in the usual Hilbert space
norm on that space.
\end{definition}

In fact, it follows from basic functional
analysis, that $L^2$-density holds automatically for
noncommutative $H^\infty$ algebras.
Indeed,
if $y \in L^2(M)$ with $y \perp A + A^*$, then
$y \in L^1(M)$, which forces $y = 0$ by the
definition of noncommutative $H^\infty$.
We will see later
that Szeg\"o's theorem implies
$L^2$-density.

\medskip

{\bf Remark.}  For a tracial subalgebra, $L^2$-density is not equivalent
to $A$ being a noncommutative $H^\infty$, or
 satisfying Szeg\"o's theorem (e.g.\
see \cite[Section 4]{HR}).

\medskip

If $A$ is a tracial subalgebra of $M$, we will write
$A_\infty$ for $[A]_2 \cap M$, where $[A]_2$ is the
norm-closure of $A$ in $L^2(M)$.

\begin{theorem} \label{hinf2}  If
$A$ is a tracial subalgebra of $M$, then
so is $A_\infty$, with canonical conditional
expectation extending that of $A$.
Also,
$A = A_\infty$, if $A$ is a noncommutative $H^\infty$
algebra.
  \end{theorem}

\begin{proof}
To see that $A_\infty$ is weak* closed in $M$,
suppose that $x \in M$ was in the
weak* closure of $A_\infty$.
If $x$ were not in $[A]_2$, then we could find
$b \in L^2(M)$ with
$b \perp [A]_2$, but $\tau(b x) \neq 0$.
But then $b \in L^1(M)$ with
$\tau(b y) = 0$ for all $y \in A_\infty$,
and consequently $\tau(b x) = 0$ since $x$
is in the weak* closure of $A_\infty$.
This contradiction shows that
 $x \in [A]_2 \cap M = A_\infty$, so that
$A_\infty$ is weak* closed.

To see that $A_\infty$ is an algebra,
one first checks that if $a \in A, b \in
A_\infty$, then $a b \in A_\infty$.
Indeed, if $(b_n) \subset A$ with $b_n \to
b$ in $L^2(M)$, then $a b_n \in A$, and
$a b_n \to ab$ in $L^2(M)$.  Thus
$a b \in [A]_2 \cap M = A_\infty$.
If $a \in A_\infty$, and if
$(a_n) \subset A$ with $a_n \to
a$ in $L^2(M)$, then $a_n b \in A_\infty$
by what we just proved, and
$a_n b \to ab$ in $L^2(M)$.  Thus again
$a b \in [A]_2 \cap M = A_\infty$.

We continue to write  $\Phi$ for the canonical
conditional expectation from $M$ to ${\mathcal D}$
extending the projection from $A$ onto ${\mathcal D}$
(e.g.\ see \cite[Theorem 5.6]{BL}), and for
the further extension to $L^p(M)$ (e.g.\ see
\cite[3.9]{MW}).   We
claim that $\Phi(ab) = \Phi(a) \Phi(b)$ for all
$a, b \in [A]_2$.   Indeed, if $a_n, b_n \in A$
with $a_n \to a$ and $b_n \to b$ in $L^2(M)$,
then $a_n b_n \to ab$ in $L^1(M)$, so that
$$\Phi(ab) \; = \; \lim_n \Phi(a_n b_n)
 \; = \; \lim_n \Phi(a_n) \Phi(b_n)
 \; = \; \Phi(a) \Phi(b) ,$$
by the continuity of  $\Phi$  on each
$L^p(M)$.   In particular, $\Phi$ is
a homomorphism on $A_\infty$.
By the argument in the second paragraph of the
proof of \cite[Proposition 2.1.4]{AIOA},
$A_\infty \cap A_\infty^* = {\mathcal D}$.

If $A$ is a noncommutative $H^\infty$, then
it is well known that $A_\infty = A$.
This follows, for example, from the fact that
noncommutative $H^\infty$ algebras are $\tau$-maximal
(see \cite[Section 2]{AIOA}).
 \end{proof}

{\bf Remark.}  It is interesting that if $A$ satisfies
 Jensen's inequality, then $A_\infty$  satisfies Jensen's formula.
To see this, note that by lines 8 and
9 of p.\ 613 of \cite{AIOA}, we need only show
that $\Delta(a) \geq \Delta(\Phi(a))$ for every $a \in A_\infty^{-1}$. To
prove this, by replacing $a$ by $a\Phi(a)^{-1}$ if necessary,
it suffices to show that $\Delta(a) \geq 1$ for every
invertible $a \in A_\infty$ with $\Phi(a) = 1$.
It is easy to approximate such $a$ in the $L^2$-norm,
and hence in $L^1$-norm, by a sequence in $A$
      with $\Phi(a_n) = 1$ for each $n$.
By \cite[III.4.10]{Tak}, $|a_n| \to |a|$ in $L^1$-norm.
If $\epsilon > 0$ is given, by Proposition \ref{Arvf} there exists
a $b \in M^{-1} \cap M_+$ with $\Delta(b) \geq 1$, such that
$$\Delta(a) + \epsilon \geq \tau(b|a|) = \lim_n \tau(b|a_n|) \geq \limsup_n
\Delta(a_n). $$
Thus, $\Delta(a) \geq \limsup_n \Delta(a_n)$, and since
Jensen's inequality holds for $A$, we have
$$ \Delta(a) \geq \limsup_n \Delta(a_n)
\geq \limsup_n \Delta(\Phi(a_n)) = \Delta(1) = 1$$
 as required.

\begin{definition} \label{pf}  We shall say that a
tracial subalgebra of $M$ has {\em partial factorization},
if whenever $b \in M_+ \cap M^{-1}$, we have
$b = |a| = |c^*|$ for some elements $a, c \in A \cap M^{-1}$, with
$\Phi(a) \Phi(a^{-1}) = \Phi(c) \Phi(c^{-1}) = 1$.
 \end{definition}

We remark that this notion is distinct from the one-sided
partial factorization considered in \cite{Pitts}.

\begin{theorem} \label{pfac}  Let $A$ be a tracial subalgebra of $M$.
\begin{itemize}
\item [(a)]  If $A$ has the unique normal state extension
property, then $A_\infty$ has partial factorization.
If, further, $A$  satisfies $L^2$-density,
then $A_\infty$ has factorization,
and is a noncommutative $H^\infty$ algebra.
\item [(b)]  If $A$ satisfies Szeg\"o's theorem, then
$A$  satisfies $L^2$-density,
and hence $A_\infty$ is a noncommutative $H^\infty$ algebra.
\end{itemize}
\end{theorem}

\begin{proof}
(a) \ Let $b \in M_+ \cap M^{-1}$ be given, and define an inner product on
$L^2(M)$ by
$$\langle f , g \rangle_b \; = \; \tau(b^{\frac{1}{2}} g^* f b^{\frac{1}{2}}) \;
, \qquad f , g \in L^2(M) .$$
The norm $\Vert \cdot \Vert_b$ induced by this inner product
is equivalent to the usual one.  Indeed,
$$\Vert f \Vert^2_b = \tau(b |f|^2) \leq \Vert b \Vert
\tau(|f|^2) \leq \Vert b \Vert \Vert b^{-1} \Vert \tau(b |f|^2)
= \Vert b \Vert \Vert b^{-1} \Vert \Vert f \Vert^2_b .$$

Let $p$ be the orthogonal projection of $1$ into the subspace $[A_0]_2$,
taken with respect to the inner product $\langle \cdot , \cdot
\rangle_b$. (Here $[A_0]_2$ is the $L^2$-closure of $A_0$ in
$L^2(M)$.)  Now  $1-p \in [A]_2$, since $p \in [A_0]_2 \subset [A]_2$.
Therefore, by an obvious approximation argument,
$a_0(1-p)  \in [A_0]_2$ for all $a_0 \in A_0$.
 Since $1-p \perp [A_0]_2$, we have
$$0 = \langle  a_0(1 - p) , 1 - p \rangle_b \; = \; \tau(b^{\frac{1}{2}}
(1 - p^*) a_0 (1 - p) b^{\frac{1}{2}}) =
\tau(a_0 (1 - p) b (1 - p^*)) .$$
By the hypothesis, combined with Lemma \ref{gcon},
we deduce that $(1 - p) b (1 - p^*) \in L^1({\mathcal D})_+$.

 Let $\epsilon > 0$ be given such that
$b \geq \epsilon 1$.  If $d \in {\mathcal D}_+$, then
$$\tau((1 - p) b (1 - p^*) d)
= \tau(d^{\frac{1}{2}} (1 - p) b (1 - p^*) d^{\frac{1}{2}})
 \geq \epsilon \tau(d^{\frac{1}{2}} (1 - p) (1 - p^*) d^{\frac{1}{2}}) .$$
By
the $L^2$-contractivity of $\Phi$ we deduce that
$$\tau((1 - p) b (1 - p^*) d) \geq \epsilon \tau(|d^{\frac{1}{2}} (1-p)|^2)
\geq \epsilon \tau(|d^{\frac{1}{2}} \Phi(1-p)|^2)
= \epsilon \tau(d) .$$
This implies that $(1 - p) b (1 - p^*)  \geq \epsilon 1$ in
$L^1({\mathcal D})$.
Thus this element has a bounded
inverse $((1 - p) b (1 - p^*))^{-1} \in {\mathcal D}$.
Set $e = ((1 - p) b (1 - p^*))^{-\frac{1}{2}} \in {\mathcal D}$,
and let $a = e(1-p) \in [A]_2$.  Since
$$1 = e (1 - p) b (1 - p^*) e = |b^{\frac{1}{2}} (1 - p^*) e|^2 ,$$
we deduce that $b^{\frac{1}{2}} (1 - p^*) e$, and consequently also
$(1 - p^*) e$ and $a = e (1-p)$, are bounded.  Since they
belong to $L^2(M)$, we deduce that $a \in M$.
Hence $a \in M \cap [A]_2 = A_\infty$.  Since $1 = a b a^*$,
and since $M$ is a finite von Neumann algebra, we also have
$1 = b a^* a$, so that $b^{-1} = |a|^2$.

Note that $\Phi(a) = e\Phi(1-p) = e$,
since $\Phi$ is $L^2$-continuous, and $p \in [A_0]_2$.
We have $a^{-1} = b(1-p^*) e$.
To see that $\Phi(a) \Phi(a^{-1}) = 1$,
first note that for any $d \in {\mathcal D}$,
$$\tau(\Phi(p b(1-p^*)) d) = \tau(p b(1-p^*) d) =
\langle d p , 1-p \rangle_b
= 0,$$
since $dp$ is in $[A_0]_2$, and $1-p \perp [A_0]_2$.
 Since this is true for all
$d \in {\mathcal D}$, we have $\Phi(p b(1-p^*)) = 0$.
Thus $$1 = e \Phi((1-p)b(1-p^*)) e
= e \Phi(a^{-1}) = \Phi(a) \Phi(a^{-1}) .$$

Notice that we now have that $e^{-1}$ is bounded, and this in
fact forces $p$ to be bounded.  Therefore $p \in M$.  However
we do not see any use for the latter fact at this point.

An analogous argument, using the inner product
$\langle f , g \rangle^b
= \tau(b^{\frac{1}{2}} f g^* b^{\frac{1}{2}})$, gives
 $b^{-1} = |c^*|^2$, for some $c \in A_\infty$, and
$1 = \Phi(c) \Phi(c^{-1})$.

Now suppose that $A + A^*$ is norm-dense in $L^2(M)$.
Then by an obvious argument,
$L^2(M) = [A]_2 \oplus [A^*_0]_2$, where
$[A^*_0]_2$ is the norm-closure of $A_0$ in  $L^2(M)$.
For any $a_0 \in A_0$, we have
$$\tau(a^{-1} a_0) = \tau(b(1-p^*) e a_0) =
\langle e a_0 , 1 - p \rangle_b = 0, $$
since $ea_0 \in [A_0]_2$, and $1-p \perp [A_0]_2$.
Thus
$$a^{-1} \in M \cap (L^2(M) \ominus [A^*_0]_2)
=  M \cap [A]_2 = A_\infty .$$
We deduce that $A_\infty$ has factorization,
and is therefore a noncommutative $H^\infty$ algebra
by Corollary \ref{logeqnon}.

(b) \ Suppose that $A$ satisfies Szeg\"o's theorem.
We will prove that $A + A^*$ is norm-dense in $L^2(M)$,
and then the result follows from (a)
and Theorem \ref{gcond}.  In fact,
if $A$ has the normal state extension property,
then we will only need Szeg\"o's theorem for $h \in M_+$, as opposed
to $h \in L^1(M)_+$.  For suppose that $k \in L^2(M)$ is
such that $\tau(k (A + A^*)) = 0$. We show that the previously
stated conditions are enough to then force $k = 0$. By the argument
in the proof of Corollary \ref{logeqnon}, we may assume that
$k = k^*$. Then $1 - k \in L^1(M)$,
so that by (\ref{Eq2}), given $\epsilon > 0$ there exists
an element $b \in M_+ \cap M^{-1}$
with $\Delta(b) \geq 1$ and $\tau(|(1 - k) b|) < \Delta(1 - k) +
\epsilon$.  By scaling, we may assume that $\Delta(b) = 1$,
and hence $\Delta(b^{-2}) = 1$ by the multiplicativity of
$\Delta$.
By Szeg\"o's theorem, there exists an $a_0 \in A_0$, and an
invertible $d_0 \in {\mathcal D}$ with $\Delta(d_0) \geq 1$,
such that
\begin{equation} \label{szr}
1 - \epsilon  = \Delta(b^{-2}) - \epsilon \geq
\tau(b^{-2} |d_0 + a_0|^2) .
\end{equation}   We will modify the argument
at the beginning of the proof, but with $b$ replaced by
$b^{-2}$.  Thus we consider the inner product
$\langle x , y \rangle_{b^{-2}} = \tau(b^{-2} y^* x)$ above,
but now we let $p$ be the projection of $d_0$ onto the
subspace $[A_0]_2$.  By a simple modification of the earlier
argument, $(d_0 - p) b^{-2} (d_0 - p)^*
\in L^1(D)_+$, and there are constants $K_0, K_1 > 0$ such that
for any $d \in {\mathcal D}_+$,
$$\tau((d_0 - p) b^{-2} (d_0 - p)^* d) \; \geq \; K_0
\tau(d^{1/2}|d_0^*|d^{1/2}) \; \geq \; K_0K_1 \tau(d).$$
We conclude, as before, that $(d_0 - p) b^{-2} (d_0 - p)^*$ has
a bounded inverse in ${\mathcal D}$, and we set
$e = ((d_0 - p) b^{-2} (d_0 - p)^*)^{-\frac{1}{2}}$.
 As before, $a = e(d_0 - p) \in A_\infty$,
$e^{-1} \in {\mathcal D}$, and $b^2 = a^* a$.   If $(a_n) \in A$ with
$a_n \to a$ in  $L^2$-norm, and if
$d \in {\mathcal D}$, then by H\"older's inequality
we have that $a_n k d
\to a k d$ in $L^1$-norm.  Thus $\tau(akd) = 0$,
and hence it follows, by a routine argument that
we have used several times already, that $\Phi(ak) = 0$.
By the $L^1$-contractivity of $\Phi$, and (\ref{Eq2}),  we have
$$\tau(|a - ak|) \geq \tau(|\Phi(a - ak)|)
=  \tau(|\Phi(a)|) = \tau(|e d_0|) \geq \Delta(e d_0)
=  \Delta(e)  \Delta(d_0) \geq  \Delta(e).$$
Hence
$$\tau(|(1 - k) b|) = \tau(|b(1 - k)|)
= \tau(||a|(1 - k)|) = \tau(|a(1 - k)|) \geq \Delta(e) .$$
On the other hand, equation (\ref{szr}) informs us that
$\Vert d_0 + a_0 \Vert^2_{b^{-2}} \leq 1 - \epsilon,$ and so
by
definition of $p$ we have
$\Vert d_0 - p \Vert^2_{b^{-2}} \leq \Vert a_0 + d_0
\Vert^2_{b^{-2}} \leq  1 - \epsilon$.  Thus,
by (\ref{Eq1}) and the definition of $\Vert \cdot \Vert_{b^{-2}}$,
 $$\Delta(e^{-2}) \leq \tau(e^{-2}) =
\Vert d_0 - p \Vert^2_{b^{-2}} \leq 1 - \epsilon .$$
By the multiplicativity of
$\Delta$, we have
$\Delta(e) \geq \frac{1}{\sqrt{1 - \epsilon}}$.
We conclude that
$$\Delta(1-k) \geq \tau(|(1 - k) b|) - \epsilon  \geq \Delta(e) - \epsilon
\geq \frac{1}{\sqrt{1 - \epsilon}} - \epsilon.$$
Since $\epsilon > 0$ is arbitrary, $\Delta(1-k) \geq 1$.
As in the proof of Corollary \ref{logeqnon}, we conclude that $k = 0$.  Thus
$A + A^*$ is norm-dense in $L^2(M)$.
\end{proof}

\begin{corollary} \label{ifjen}   If $A$ is a
$\tau$-maximal tracial subalgebra of $M$
satisfying
the unique normal state extension property,
then $A$
is a noncommutative $H^\infty$ algebra.
\end{corollary}

\begin{proof}    Let $A$ be a
$\tau$-maximal tracial subalgebra of $M$.
If $b \in A_\infty$, then there exists a
sequence $(a_n) \subset A$ with
$L^2$-limit $b$.  If $c \in A_0$, then
$a_n c \to b c$ in $L^2(M)$, and
$$\tau(b c)  \; = \;  \lim_n \tau(a_n c) \; = \; 0 .$$
Thus $b \in A$.  Therefore $A_\infty = A$.
By the last part of the proof of Theorem  \ref{pfac} (a),
$\tau(a^{-1} c) = 0$ too, so that $a^{-1} \in A$.
Thus $A$ has factorization, and is therefore a
noncommutative $H^\infty$ algebra
by Corollary \ref{logeqnon}.
\end{proof}

\section{Maximal algebras}

If $A$ is a tracial subalgebra of $M$, with canonical
projection $\Phi : A \to {\mathcal D}$, then as
observed in \cite{BL},
$\Phi$ extends canonically to a conditional expectation
from $M$ onto ${\mathcal D}$, which we continue to write as $\Phi$,
such that $\tau \circ \Phi = \Phi$.  We say that $A$ is
{\em maximal} if there is no properly larger tracial
subalgebra of $M$ with conditional expectation
$\Phi$.  Equivalently, there is no properly larger tracial
subalgebra $B$ of $M$ whose
conditional  expectation onto $B \cap B^*$ extends
the one on $A$.
In this section, we will also consider tracial subalgebras
which satisfy $L^2$-density
(resp.\ the unique normal state extension property (\ref{cond})).
  In this case, any larger tracial subalgebra clearly also has
these properties, and so  $A$ is maximal if and only
if it is maximal among the tracial subalgebras
with these properties and with conditional expectation
$\Phi$.

Exel showed that noncommutative $H^\infty$ algebras
are automatically maximal \cite{E}.  We will
have to establish this for algebras satisfying some of our
other conditions studied above.

\begin{lemma} \label{max}  Let $A$ be a
tracial subalgebra of $M$ with conditional expectation
$\Phi$, which satisfies
$L^2$-density.
Then there is a unique
largest tracial subalgebra of $M$ with conditional expectation
$\Phi$ containing $A$, namely:
$$\{ x \in M \, : \, \tau(x a) = 0 \; \text{for all}
\; a \in A_0 \} . $$
In particular,  if
$A$ is maximal and satisfies $L^2$-density,
then $A$ is
$\tau$-maximal, that is $A
 = \{ x \in M \, : \, \tau(x a) = 0 \; \text{for all}
\; a \in A_0 \} . $
\end{lemma}

\begin{proof} This follows by adapting the arguments in \cite[Section
2.2]{AIOA}.  In particular,
using the $L^2$-density as opposed to weak* density
to simplify the argument from about the middle of p.\ 583 onwards.
Note that $\tau(x a) = 0$ for all $a \in A_0$ if and only if
$\tau(x a d) = \tau(\Phi(x a)d) = 0$ for all $a \in A_0$,
and if and only if $\Phi(x a) = 0$ for all $a \in A_0$.
 \end{proof}

\begin{theorem} \label{max2}  Let $A$ be a
tracial subalgebra of $M$ which satisfies
Szeg\"o's theorem (or simply
$L^2$-density and the normal state extension property).
Then $A$ is maximal,
and $A$ is a noncommutative $H^\infty$ algebra.
\end{theorem}

\begin{proof}
It suffices, by Lemma \ref{max}
and Corollary \ref{ifjen}, to show that $A$ is maximal.
To do this we follow the proof of the
main theorem in \cite{E}, indicating how it may be adapted to show that
$A$ agrees with the maximal algebra guaranteed by Lemma \ref{max}.
Specifically, we note that if $\xi \in L^2(M)$, and if
$\langle  x \xi , \xi \rangle = \langle \Phi(x) \xi , \xi \rangle$
for all $x \in A + A^*$, then
$\tau(x \xi \xi^*) = \langle  x \xi , \xi \rangle =
0$ if $x \in A_0$.  By (\ref{cond}), we deduce that
$\xi \xi^* \in L^1({\mathcal D})$.  It follows that for any
$x \in M$, we have
$$\langle \Phi(x) \xi , \xi \rangle =  \tau(\Phi(x) \xi \xi^*) =
\tau(\Phi(x \xi \xi^*)) = \tau(x \xi \xi^*) = \langle x \xi , \xi
\rangle,$$
as needed for the proof in the middle of page 779 of \cite{E}
to proceed.  The next obstacle one encounters
is that the first centered equation on \cite[p.\ 780]{E} is
not clear for the vector $\xi_1$ there.  However,
because $A$ satisfies $L^2$-density,
this equation is easily seen to be equivalent to the
next group of centered equations.
In the notation of that paper, we have
$u \overline{{\mathcal D}
 \xi_1} \subset \overline{{\mathcal D} \delta} = L^2({\mathcal D})$.
Since $u$ commutes with ${\mathcal D}$, it is easy to see
that $u \xi_1$ is a separating vector for the action of ${\mathcal D}$ on
$L^2({\mathcal D})$.   It follows, by \cite[Exercise 9.6.2]{KR}
 for example,
that $u \xi_1$ is a cyclic vector for the ${\mathcal D}$ action.
Thus $\overline{{\mathcal D} u \xi_1} = \overline{{\mathcal D} \delta}$.
That is, $\delta = \lim_n \, d_n u \xi_1$ for a sequence $(d_n) \subset
{\mathcal D}$.  We have
$$u \overline{A_0 \xi_1} \subset \overline{A_0 u \xi_1}
\subset \overline{A_0{\mathcal D} \delta} \subset \overline{A_0 \delta}
.$$
Conversely,
$$a_0 \delta \, = \, \lim_n \, a_0 d_n u \xi_1 \; \in \; u \overline{A_0
\xi_1} \; , \qquad a_0 \in A_0 .$$
Thus $u \overline{A_0 \xi_1} = \overline{A_0 \delta}$.  Similarly,
$u \overline{A^*_0 \xi_1} = \overline{A^*_0 \delta}$.
 As we said above, this gives the first centered equation on \cite[p.\
780]{E}.
The rest of the proof is unchanged.
   \end{proof}

\section{Other characterizations of noncommutative
$H^\infty$ algebras}

We now turn to the remaining items in the list in \cite[Section 3]{SW}
of conditions equivalent to logmodularity.  We recall
that a simply (right) invariant subspace of $L^2(M)$,
is a closed subspace $M$ of $L^2(M)$ such that
$M A \subset M$, and the closure of
$M A_0$ is properly contained in $M$.   We have:

\begin{theorem} \label{rest}  Let $A$ be a tracial subalgebra of $M$, and
consider the following properties:
\begin{itemize} \item [(a)] (Invariant subspace theorem) \
 Every simply right invariant subspace of $L^2(M)$
is of the
form $u [A]_2$, for a unitary $u$ in $M$,
\item [(b)]  (Beurling-Nevanlinna factorization) \
Whenever  $f \in L^2(M)$ with $f \notin [f A_0]_2$,
then $f = uh$, for a unitary $u$ in $M \cap [f A]_2$ and an $h$ with
$[h A]_2 = [A]_2$,
\item [(c)]   $A_\infty$ is a noncommutative
$H^\infty$ algebra,
 \item [(d)]   $A$ is a noncommutative
$H^\infty$ algebra.
\end{itemize}
Then {\rm (a)} $\Rightarrow$ {\rm (b)} $\Rightarrow$
{\rm (c)}.  If $A$ is antisymmetric
(that is, ${\mathcal D}$ is one-dimensional), then
{\rm (d)} $\Rightarrow$ {\rm (a)}.
If $A$ has the unique normal state extension
property, then {\rm (c)} $\Rightarrow$ {\rm (d)}.
\end{theorem}

\begin{proof}  That {\rm (a)} $\Rightarrow$ {\rm (b)} $\Rightarrow$
{\rm (c)} follows just as in \cite{SW}, for example, with
insignificant modifications.   One also needs to use the fact about
invertibility in a finite algebra stated at the end of
Section 1.   The assertion about
antisymmetric algebras is also essentially just as in
\cite{SW} (see also \cite{Kam}).    Finally,
if $A_\infty$ is a noncommutative
$H^\infty$ algebra then it satisfies $L^2$-density.
Therefore it is clear that $A$ satisfies $L^2$-density.
Appealing to Theorem \ref{max2}, we obtain {\rm (d)}.
 \end{proof}

As the last Theorem shows, it is of interest to know whether
$A_\infty$ being a  noncommutative
$H^\infty$ algebra implies that $A$ is a noncommutative
$H^\infty$ algebra.   We end the paper with another
 sufficient condition under which this holds.

\begin{lemma}
Let $A$ be a tracial subalgebra of $M$. Then
\begin{center}
$\{ x \in L^1(M) \, : \, \tau(x a) = 0 \; \text{for all}
\; a \in A \} =$
\end{center}
\begin{center}
$\{ x \in L^1(M) \, : \, \tau(x a) = 0 \; \text{for all} \; a \in
A_0 \} \cap \mathrm{Ker}(\Phi).$
\end{center}
(Here we
identify $\Phi$ with its extension to $L^1(M)$.) Moreover for any
$y \in \{ x \in L^1(M) \, : \, \tau(x a) = 0 \; \text{for all} \;
a \in A_0 \}$, we have $\Phi(by) = \Phi(b)\Phi(y)$ for all $b \in A$.
\end{lemma}

\begin{proof} Let $x \in L^1(M)$ be given. It is easy to see that
$\tau(ax) = 0$ for all $a \in A$ if and only if $\tau(ax) = 0$ for
all $a \in A_0$ and $\tau(dx) = 0$ for all $d \in {\mathcal D}.$ Since
for any $d \in {\mathcal D}$ we have $\tau(dx) = \tau(\Phi(dx)) =
\tau(d\Phi(x))$ and since $\Phi(x) \in L^1({\mathcal D})$, the first
claim follows.

To see the second claim, let $y \in \{ x \in L^1(M) \, : \,
\tau(x a) = 0 \; \text{for all} \; a \in A_0 \}$ be given. For any
$d \in {\mathcal D}$ and $a \in A_0$, we then have $$\tau(d\Phi(ay)) =
 \tau(\Phi(day)) = \tau((da)y) = 0.$$ Since $\Phi(ay) \in L^1({\mathcal D})$,
this suffices to show that $\Phi(ay) = 0$ for all $a \in A_0$. But then
given any $b \in A$, we will surely have $$\Phi(by) = \Phi(\Phi(b)y) +
\Phi(b - \Phi(b))y) = \Phi(\Phi(b)y) = \Phi(b)\Phi(y).$$
This is the desired identity.
\end{proof}

By analogy with the commutative
context, the space $\widetilde{M}$ mentioned at the end of the
introduction may also be equipped with a topology of
\emph{convergence in measure}, in such a way that each $L^p(M,\tau)$
injects continuously into $\widetilde{M}$ (see
\cite{Terp,FK,Tak2,Nel} for details).
With respect
to this topology, $\widetilde{M}$
becomes a complete Hausdorff topological
$\ast$-algebra with respect to the `strong' sum and product.

\begin{corollary}
Let the canonical extension of $\Phi$ to $L^1(M)$ be continuous with
respect to the topology of convergence in measure. Then
\begin{center}
$\{ x \in L^1(M) \, : \, \tau(x a) = 0 \; \text{for all}
\; a \in A_0 \} =$
\end{center}
\begin{center}
$\{ x \in L^1(M) \, : \, \tau(x a) = 0 \; \text{for all}
\; a \in (A_\infty)_0 \}.$
\end{center}
\end{corollary}

\begin{proof} Since $A_0 \subset (A_\infty)_0$, the one inclusion is
  trivial.
 For the converse suppose that we are given $y \in \{ x \in
L^1(M) \, : \, \tau(x a) = 0 \; \text{for all} \; a \in A_0 \}$. By
the lemma we then have that $\Phi(by) = \Phi(b)\Phi(y)$ for all $b \in
A$.
  Given any $a_0 \in (A_\infty)_0$ we may now select $(a_n) \subset
  A_0$ so
  that $(a_n)$ converges to $a_0$ in $L^2$-norm.  Thus
 $(a_n)$ converges to $a_0$ in the topology of convergence
in measure (e.g.\ see \cite[Theorem 5]{Nel}).
 Since $\widetilde{M}$ is a topological algebra
in this topology \cite{Nel,Terp}, $a_ny \rightarrow a_0 y$ in this topology.
Therefore if the extension of $\Phi$ to $L^1(M)$ is indeed continuous
with respect to this topology, then
$$0 = \lim_n \Phi(a_n)\Phi(y) =
 \lim_n \Phi(a_ny) = \Phi(a_0y) . $$
 This clearly forces $0 = \tau(\Phi(a_0y))
= \tau(a_0y)$ as required. \end{proof}

 \begin{theorem}
Suppose that the canonical extension of $\Phi$ to $L^1(M)$ is
continuous with
respect to the topology of convergence in measure.
Then $A_\infty$ is a
noncommutative $H^\infty$ algebra
if and only if any one of the equivalent conditions
{\rm (a)}--{\rm (f)} in Theorem {\rm \ref{main}} holds.
In particular, if $A_\infty$ is a
noncommutative $H^\infty$ algebra, then $A = A_\infty$.
\end{theorem}

\begin{proof}
Suppose that $A_\infty + A_\infty^*$ is weak* dense in $M$, and that the
extension of $\Phi$ to $L^1(M)$ is continuous in the
topology of convergence in measure. Let $g \in L^1(M)$ be
given with $g \perp A + A^*$. To prove the result,
by Theorem \ref{hinf2}, it is enough to show that
then $g = 0$.   It clearly suffices to show that  if
$g \perp A$ then $g \perp A_\infty$.  This is in turn a trivial
consequence of the preceding results applied to $g^*$.
\end{proof}

{\bf Closing Remark.}    Although most results in
Arveson's paper \cite{AIOA} are stated for
{\em finite subdiagonal subalgebras}  of
von Neumann algebras with a faithful normal tracial state,
he also considers {\em subdiagonal subalgebras} of
general von Neumann algebras.  It would be interesting
if there was some way to extend some of our results to this context.
See e.g.\ \cite{Xu} for recent work on the question of
maximality for this larger class of algebras.

\end{document}